\documentclass{amsart}
\theoremstyle{plain}
\usepackage{amsfonts}
\usepackage{graphicx}
\usepackage{enumitem}
\usepackage{amssymb}
\def\str#1{\mathbf {#1}}
\def\todo#1{}
\newtheorem{theorem}{Theorem}[section]

\theoremstyle{definition}
\newtheorem{definition}[theorem]{Definition}

\newtheorem*{remark*}{Remark}
\newtheorem*{claim*}{Claim}
\theoremstyle{remark}

\def\str#1{\mathbf {#1}}
\def\rel#1#2{R_{\mathbf{#1}}^{#2}}
\begin{document}


\title{Big Ramsey degrees of 3-uniform hypergraphs}


\author{M. Balko}
\address{Department of Applied Mathematics (KAM), Charles University, Ma\-lo\-stransk\'e n\'am\v est\'i 25, Praha 1, Czech Republic}
\email{balko@kam.mff.cuni.cz}

\author{D. Chodounsk\'y}
\address{Department of Algebra, Charles University, Ma\-lo\-stransk\'e n\'am\v est\'i 25, Praha 1, Czech Republic}
\email{chodounsky@math.cas.cz}

\author{J. Hubi\v cka}
\address{Department of Applied Mathematics (KAM), Charles University, Ma\-lo\-stransk\'e n\'am\v est\'i 25, Praha 1, Czech Republic}
\curraddr{}
\email{hubicka@kam.mff.cuni.cz}
\urladdr{}

\author{M. Kone\v{c}n\'{y}}
\address{Department of Applied Mathematics (KAM), Charles University, Ma\-lo\-stransk\'e n\'am\v est\'i 25, Praha 1, Czech Republic}
\email{matej@kam.mff.cuni.cz}

\author{L. Vena}
\address{Department of Applied Mathematics (KAM), Charles University, Ma\-lo\-stransk\'e n\'am\v est\'i 25, Praha 1, Czech Republic}
\email{lluis.vena@gmail.com}



\thanks {Supported  by  project  18-13685Y  of  the  Czech  Science Foundation (GA\v CR) and by Charles University projects Progres Q48 and GA UK No 378119.}



{}

\begin{abstract}
Given a countably infinite hypergraph $\mathcal R$ and a finite hypergraph $\mathcal A$, the \emph{big Ramsey degree} of $\mathcal A$ in $\mathcal R$ is
the least number $L$ such that, for every finite $k$ and every $k$-colouring of the embeddings of $\mathcal A$ to $\mathcal R$, there exists an embedding $f$ from $\mathcal R$ to $\mathcal R$ such that all the embeddings of $\mathcal A$ to the image $f(\mathcal R)$ have at most $L$ different colours.

We describe the big Ramsey degrees of the random countably infinite 3-uniform hypergraph, thereby solving a question of Sauer. We also give a new presentation of the results of Devlin and Sauer on, respectively, big Ramsey degrees of the order of the rationals and the countably infinite random graph.
Our techniques generalise (in a natural way) to relational structures and give new examples of Ramsey structures (a concept recently introduced by Zucker with applications to topological dynamics).
\end{abstract}

\maketitle

\section{Introduction}
We consider graphs, hypergraphs and orders as special cases of (relational)
structures defined more formally below.
Given structures $\str{A}$ and $\str{B}$, we denote by $\str{B}\choose \str{A}$ the set
of all embeddings from $\str{A}$ to $\str{B}$. We write $\str{C}\longrightarrow (\str{B})^\str{A}_{k,l}$ to denote the following statement:
for every colouring $\chi$ of $\str{C}\choose\str{A}$ with $k$ colours, there exists
an embedding $f:\str{B}\to\str{C}$ such that $\chi$ does not take more than $l$ values on $f(\str{B})\choose \str{A}$.
Given a class $\mathcal K$ of structures, the \emph{(small) Ramsey degree} of $\str{A}$ in $\mathcal K$
is the least $l\in \mathbb N\cup\{\infty\}$ such that for every $\str{B}\in \mathcal K$ and $k\in \mathbb N$ there exists
$\str{C}\in \mathcal K$ such that $\str{C}\longrightarrow (\str{B})^\str{A}_{k,l}$.
For a countably infinite structure $\str{B}$ and its finite substructure $\str{A}$, the \emph{big Ramsey degree} of $\str{A}$ in $\str{B}$ is
the least number $L\in \mathbb N\cup \{\infty\}$ such that $\str{B}\longrightarrow (\str{B})^\str{A}_{k,L}$ for every $k\in \mathbb N$.

The class $\mathcal K$ is said to be a \emph{Ramsey class} if the Ramsey degree of every $\str{A}\in \mathcal K$ is 1.
Ramsey classes are the main topic of interest
of the structural Ramsey theory. To this date, there are many Ramsey classes known. Examples relevant for our presentation include the class of all finite linear
orders, the class of all finite graphs equipped with a linear order on their vertices, and
the class of all finite 3-uniform hypergraphs with a linear order on the vertices.  The first example is a consequence of Ramsey's theorem. Other two examples are consequences of the Ne\v set\v ril--R\"odl theorem; see~\cite{Hubicka2016}.

\medskip

Contrary to the study of Ramsey classes, which has been a very active area recently,
big Ramsey degrees are much less understood and very few non-trivial examples are known.
One of the difficulties is the fact that Ramsey degrees are often determined by the number of non-isomorphic orderings (thus being 1 for classes of structures with ordering on vertices), while
big Ramsey degrees are surprisingly rich (as discussed later). 
The interest in the area was recently renewed by the work of Zucker~\cite{Zucker2017}, which gives a good equivalent of Ramsey expansion for big Ramsey degrees, and
by deep results of Dobrinen~\cite{dobrinen2019ramsey3, dobrinen2019ramsey} about big Ramsey degrees of Henson graphs (see~\cite{dobrinen2019ramsey}).

We give a new presentation of two classical results in the area---the big Ramsey degrees of the (natural) order of the rationals by Devlin~\cite{devlin1979} and of the Rado graph by Sauer~\cite{Sauer2006}. Solving a question of Sauer\footnote{Personal communication, 2014}, we generalise this to the description of big Ramsey degrees of the random 3-uniform hypergraph.
This is a contribution to the ongoing efforts to give a more systematic treatment to this, so far, very mysterious area.

\section{Preliminaries}
We use the standard model-theoretic notion of relational structures.
Let $L$ be a language with relational symbols $\rel{}{}\in L$ each having its {\em arity}.
An \emph{$L$-structure} $\str{A}$ on $A$ is a structure with {\em vertex set} $A$, relations $\rel{A}{}\subseteq A^r$ for every symbol $\rel{}{}\in L$ of arity $r$.  If the set $A$ is finite we call $\str A$ a \emph{finite structure}. We consider only structures with finitely many or countably infinitely many vertices.

We discuss several special cases of relational structures, in particular:
\begin{enumerate}
\item graphs (where the language $L$ consists of one binary relation $E$, the adjacency relation), 
\item orders (where $L$ consists of one binary relation denoted by $\leq$, $\preceq$ or $\sqsubseteq$),
\item 3-uniform hypergraphs (where $L$ consists of one ternary relation $\mathcal E$),
\item and combination of these (where one vertex set is equipped with, say, a graph and an order).
\end{enumerate}
Since we work with structures on multiple languages, we will list the vertex
set along with the relations of the structure, i.e., $(P,\leq)$ for partial orders, $G=(V,E)$ for
graphs and $\mathcal H=(V,\mathcal E)$ for hypergraphs, $(V,E,\leq)$ for a graph with order on the vertices, etc.

Given two $L$-structures $\str{A}$ and $\str{B}$, a function $f:A\to B$ is an \emph{embedding} $f:\str{A}\to\str{B}$ if it is
injective and for every $\rel{}{}\in L$ of arity $r$, $$(v_1,v_2,\ldots, v_r)\in \rel{A}{}\iff (f(v_1),f(v_2),\ldots, f(v_r))\in \rel{B}{}.$$ 
As usual in the structural Ramsey theory, given an embedding $f:\str{A}\to\str{B}$ we will call the image of $\str{A}$ in $\str{B}$ a \emph{copy} of $\str{A}$ in $\str{B}$.
We say that $\str{A}$ and $\str{B}$ are \emph{isomorphic} if there is an embedding $f:\str{A}\to \str{B}$ that is onto.

\medskip
A (rooted set-theoretic) \emph{tree} is a partially ordered set $T=(V,\sqsubseteq)$ with
a minimal element, called \emph{root}, such that for every $a\in V$ the down-set
$\{b:b\sqsubseteq a\}$ is well-ordered by $\sqsubseteq$.
In all the trees discussed here, for every pair of vertices $a,b\in V$, there is a unique vertex $a\wedge b$ called their \emph{meet} (or \emph{nearest common ancestor}).

Given $S\subseteq V$, the \emph{subtree of $T$ generated by $S$}
is $T'=(S',\sqsubseteq\restriction_{S'})$, where $S'$ is the minimal subset of $V$
closed under the meet operation such that $S\subseteq S'$. By
$\sqsubseteq\restriction_{S'}$ we denote the partial order $\sqsubseteq$ restricted to $S'$. If $S=S'$ we simply say that $S$ is a \emph{subtree} of $T$.
Note that this differs from the graph-theoretic subtree, as some vertices in $T$ between two comparable elements in $S'$ might not belong to $S'$.

Given a linear order $(V,\leq)$ and a well-preorder $(V,\preceq)$ (think of a well-order first), we denote by
$T(V,\leq,\preceq)$ a tree $(V,\sqsubseteq)$ defined by putting $a\sqsubseteq b$ if and only if
$a\preceq b$ and there is no $c$ with $a < c < b$ and $c\preceq a,b$.
If $(V,\sqsubseteq)$ is not a tree, then we leave $T(V,\leq,\preceq)$ undefined.
In computer science, this corresponds to the binary search tree for $(V,\leq)$ when items are inserted in the order given by $\preceq$.  It is easy to
see that if $(V,\leq)$ is a linear order and  $(V,\preceq)$ is a well-order, then
$T(V,\leq,\preceq)$ is always defined.

 We denote $[n]=\{1,2,\ldots,n\}$. 

\section{Big Ramsey degrees of the order of rationals}
\label{sec:1}
The result of Devlin~\cite{devlin1979} on the big Ramsey degrees of the order of rationals is shown using an application of the
Milliken tree theorem which is an infinitary Ramsey statement about infinite
finitely branching trees (see e.g.~\cite{dodos2016}).  The vertices of this
tree are described by finite $\{0,1\}$-sequences which corresponds to the rationals in the
interval $[0,2)$ in the natural way. Big Ramsey degrees are then given by sets
of such $\{0,1\}$-sequences satisfying additional axioms.  This approach was later
generalised to graphs~\cite{Sauer2006}. Generalising this approach further
however leads to numerous technical difficulties.  For this reason we do not describe details of this approach and rather discuss an alternative interpretation of these results which is related to the approach used earlier by
Sauer~\cite{Sauer98} to describe edge partitions in the triangle-free graph.

It may come as a surprise that the big Ramsey degree of $n$ tuples in the
order $(\mathbb Q,\leq)$ ir more than one for every $n>1$, while the small Ramsey degrees of linear
orders is one (as a direct consequence of Ramsey theorem). This can however be shown by the following
procedure which colours $n$-tuples in $(\mathbb Q,\leq)$ in a way that every
copy of $(\mathbb Q,\leq)$ in $(\mathbb Q,\leq)$ contains $n$-tuples with
many different colours (and, in fact, maximising the number of colours as shown by Theorem~\ref{thm:devlin}).
This construction is in a way an essential element in understanding the nature of big Ramsey degrees and, as we show later, it can be generalised
to graphs and hypergraphs, too.

\medskip

Fix an
enumeration $x_0,x_1,\ldots$. This gives a well-order $(\mathbb Q,\preceq)$ by letting $x_i\preceq x_j$ whenever $i\leq j$.
We then use $T^\preceq_\mathbb Q$ to denote the tree $T(\mathbb Q,\leq,\preceq)$. It follows from the density of $(\mathbb Q,\leq)$
that $T^\preceq_\mathbb Q$ is an infinite binary tree (regardless of the choice of the enumeration).

\begin{figure}
\centering
\includegraphics{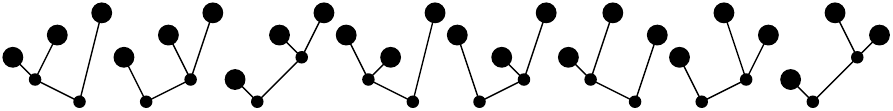}
\caption{An example of structures determining big Ramsey degrees of a linear order in $(\mathbb Q,\leq)$.}
\label{fig:1}
\end{figure}

Given $S=\{s_1, \dots, s_n\}\subset \mathbb Q$ with $s_1< \cdots < s_n$, consider the subtree $T_S$
of $T^\preceq_\mathbb Q$ generated by $S$. Since
$T^\preceq_\mathbb Q$ is a binary tree, it follows that $T_S=(S',\sqsubseteq)$ has at most $n-1$
additional vertices in the meet closure of $S$ and is thus finite.
A \emph{shape} of $S$ is the isomorphism type of the structure $(S',\leq\restriction_{S'},\preceq\restriction_{S'})$ (see Figure~\ref{fig:1} for examples of shapes of subtrees generated by 3-tuples).

Fix $n \in \mathbb{N}$ and consider the finite colouring of $n$-tuples in $(\mathbb Q,\leq)$ according to their shape (thus every shape has a unique colour assigned).
What is the least number of colours a copy of $(\mathbb Q,\leq)$ in $(\mathbb Q,\leq)$ can have?
One can observe that certain shapes can be avoided in the copy. For example subsets
  $S$ of $\mathbb Q$ for which there exist $a\neq b\in S$ such that $a\sqsubset b$. However
for $n>1$ every copy of $(\mathbb Q,\leq)$ in $(\mathbb Q,\leq)$ contains $n$-tuples of multiple colours (or shapes). 
This can be easily verified for $n=2$.
Because $\leq$ is dense and $\preceq$ is well-order it is not difficult  to show that for $2$-tuples there are 2 different shapes: every copy of $(\mathbb Q,\leq)$ in $(\mathbb Q,\leq)$ will include pairs of vertices $a<b$ where $a\wedge b\preceq a\preceq b$ and where $a\wedge b\preceq b\preceq a$.  Observe also that this argument uses the fact that $\mathbb Q$ is infinite and this is why small Ramsey degrees behave differently.

One strategy to reduce  the number of different shapes let us consider a copy
$(\mathbb{Q}',\leq)$ of $(\mathbb Q,\leq)$ in $(\mathbb Q,\leq)$ with the property that there are no $a\sqsubset b\in \mathbb Q'$, i.e., it is an \emph{antichain} in $T^\preceq_\mathbb Q$ (such a copy is not hard to construct).
Consider $S\subseteq \mathbb Q'$ with $|S|=n$ and the subtree $T_S=(S',\sqsubseteq\restriction_{S'})$
of $T^\preceq_\mathbb Q$ generated by $S$. Clearly, all vertices of $S$ are
leafs of $T_S$ and there are $n-1$ additional meet vertices.
On the other hand, and using a similar argument as the one used to show that there are two different shapes of 2-tuples, one can convince oneself that any embedding of $(\mathbb Q,\leq)$ to itself (thus preserving the ordering) contains at least one instance of each shape 
$(S',\leq\restriction_{S'},\preceq\restriction_{S'})$, such that
$|S'|=2n-1$ (so the $n$ elements of $S$ are leafs in the subtree $T_S$ of $T_{\mathbb Q}^{\preceq }$), and hence the big Ramsey degree is at least the number of such shapes.

The importance of this interaction between the linear order $\leq$ with the well-order $\preceq$ motivates the following definition.
\begin{definition}
\label{defn:1}
We say that $(X,\leq)$ and $(X,\preceq)$ is a \emph{compatible pair of orders} of $X$ if and only if $\leq$ is a linear order, $\preceq$ is a preorder, the tree $T(X,\leq,\preceq)$ is defined, and every vertex $a\in X$ is either a leaf or has 2 sons. 
\end{definition}
We note that $\preceq$ will always be a linear well-order in this section.
We will, however, need preorders later.
With the preceding argument regarding big Ramsey degrees, and
equipped with the definition of compatible orders,
 Devlin{'}s result can be
reformulated as follows.

\begin{theorem}[Devlin 1979~\cite{devlin1979}]\label{thm:devlin}
The big Ramsey degree of a suborder of $(\mathbb Q,\leq)$ with $n$ vertices
is precisely the number of non-isomorphic structures $([2n-1],\leq,\prec)$
such that $(\leq,\preceq)$ is a compatible pair of linear orders of $[2n-1]$.
\end{theorem}

Additionally, Devlin showed that the number of big Ramsey expansions of $n$-tuples
is $\tan^{(2n-1)}(0)$, the $(2n-1)$st derivative of the tangent evaluated at $0$.
The big Ramsey degrees of $n$-tuples are thus $1, 2, 16, 272, 7936, 353792, 22368256$ for $n=1,\dots,7$, respectively.
Eight out of the 16 types of linear orders with 3 vertices are depicted in Figure~\ref{fig:1}, the other 8 can be obtained by flipping the figure along the $y$ axis. We will always draw $\leq$ form left to right and $\preceq$ from bottom to top.

\section{Big Ramsey degrees of the Rado graph}
\label{sec:rado}
The \emph{Rado graph} (or the \emph{countable random graph}) is the, up to isomorphism, unique countable graph $R$ with the \emph{extension property}:
for any two finite disjoint sets $U$ and $V$ of vertices of $R$ there exists a vertex connected to all $a\in U$ and to no $b\in V$; see, for example,~\cite{Sauer2006}. 

To understand the big Ramsey degrees of finite graphs in $R$ 
we can again give a colouring of $R$ which maximises the number of colours each copy must have.
 We start by enumerating its vertices 
and defining a well-order $\preceq$.  Given a finite set of vertices $S=\{s_1,\ldots, s_n\}$ and a vertex $a$ not in $S$, the \emph{$S$-type} of $a$
is the set of all vertices $b\notin S$ of $R$ such that for every $1\leq i\leq n$ vertex $a$ is connected to $s_i$ if and only if $b$ is connected to $s_i$.
For every choice of $S$ there are $2^n$ different $S$-types of vertices of $R$ not in $S$ (given by their adjacency to vertices in $S$).

Now we can construct the \emph{tree of types} $T^\preceq_R$ of $R$ where
\begin{enumerate}
\item the vertices of $T^\preceq_R$ are all the $S$-types for every choice of $S$ being a finite initial segment of $\preceq$,
\item for two types $U,V$ we put $U\sqsubset V\iff V\subset U$.
\end{enumerate}
This is an infinite binary branching tree rooted in the $\emptyset$-type where every vertex $u$ of $R$ is associated with a unique vertex $U$ of $T^\preceq_R$ ($U$ is the inclusion minimal type in $T^\preceq_R$ containing $u$). Notice however that not every vertex of $T^\preceq_R$ has a vertex of $R$ associated this way.

By fixing an arbitrary order of the immediate successors of every vertex of the tree
this in turn leads to of a dense order $\leq$ of vertices of $R$
and thus it may not be a complete surprise that the big Ramsey degrees
of substructures of $R$ are, in fact, refinements of those of $(\mathbb Q,\leq)$. Nevertheless, it took 27 years to prove this fact.

We can refine results explained in Section~\ref{sec:1} as follows. 

\begin{figure}
\centering
\includegraphics{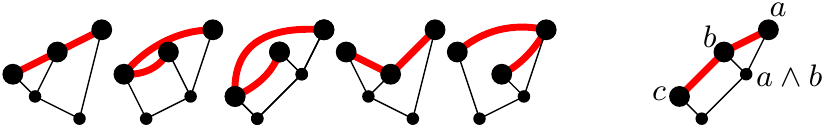}
\caption{Examples of structures determining the big Ramsey degrees of a graph path (left) and a non-example (right): $a\preceq (b\wedge c)$ but $\{a,c\}\notin E$, $\{b,c\}\in E$.}
\label{fig:2}
\end{figure}

\begin{definition}
Let $(\leq,\preceq)$ be a pair of compatible orders of a set $V'$, let $V$ be the set of leaf vertices of $T(V',\leq,\preceq)$, and let $G=(V,E)$ be a graph.
We say that $G$ is \emph{compatible} with $T(V',\leq,\preceq)$ if, for every triple $a,b,c\in V$ of distinct vertices
 satisfying $c\preceq (a\wedge b)$, we have $\{a,c\}\in E$ if and only if $\{b,c\}\in E$. 
\end{definition}
We thus annotate the tree introduced in Definition~\ref{defn:1} by a graph on its leaves in a way that it represents a subtree of $T^\preceq_R$.  See Figure~\ref{fig:2}.
Sauer's result can now be stated in an analogy to Theorem~\ref{thm:devlin}.

\begin{theorem}[Sauer 2006~\cite{Sauer2006}]\label{thm:sauer}
The big Ramsey degree of a graph $G=([n],\allowbreak E)$ in $R$
is determined by the number of non-isomorphic structures
$([2n-1],E,\leq,\preceq)$ where $(\leq,\preceq)$ is a pair of compatible linear orders
of $[2n-1]$ and $G$ is compatible with $T([2n-1],\leq,\preceq)$.
\end{theorem}
\section{Big Ramsey degrees of the $3$-uniform hypergraph}
Similarly to the Rado graph, the random countable 3-uniform hypergraph $\mathcal R$ is defined as the up to isomorphism unique 3-uniform hypergraph with the
\emph{(hypergraph) extension property}: for every pair of finite disjoint sets $U$ and $V$ of two-element subsets of $\mathcal R$, there is vertex $a$ such that $\{a,b,c\}$ is a
hyper-edge for every $\{b,c\}\in U$ and for no $\{b,c\}\in V$. 

We adapt the approach developed in Section~\ref{sec:rado}.  Enumerate the vertices of $\mathcal R$ and obtain a well-order $\preceq$.
Given a finite set of vertices $S=\{s_1,s_2,\ldots,s_n\}$ the \emph{$S$-type} of a vertex $a$ is the set of all vertices $b\notin S$ of $\mathcal R$ such that for every $1\leq i<j\leq n$ either both $\{a,s_i,s_j\}$ and $\{b,s_i,s_j\}$ are hyperedges of $\mathcal R$ or none of them is.
This let us to construct the tree of types $T^\preceq_\mathcal R$ in a complete analogy as we obtained the tree of types of the Rado graph.

This time, however, the branching degree of types
increases with $\preceq$: there is one $\emptyset$-type, one $\{s_1\}$-type (both consisting of all vertices of $\mathcal R$), two $\{s_1,s_2\}$-types and each of them refines to four $\{s_1,s_2,s_3\}$-types each of which are further refined to eight $\{s_1,s_2,s_3,s_4\}$-types and so on.

\medskip

This is not the only difference. For each vertex $a$ of $\mathcal R$ one can
consider a \emph{neighbourhood graph} $R_a$ defined on vertices of $\mathcal R$ without $a$ such
that $\{b,c\}$ forms an edge of $R_a$ if and only if $\{a,b,c\}$ is a hyper-edge of $\mathcal
R$. $R_a$ is isomorphic to the Rado graph and, since
big Ramsey structures preserved by fixing a vertex, the big Ramsey degrees in $\mathcal R$ must
account for those of $R_a$ for each choice of $a$. The graph $R_a$ has its own tree of types $T^\preceq_{R_a}$ constructed by the method described in Section~\ref{sec:rado} and introduces additional structure which we need to account for while determining the big Ramsey degrees of $\mathcal R$.

\begin{figure}
\centering
\includegraphics{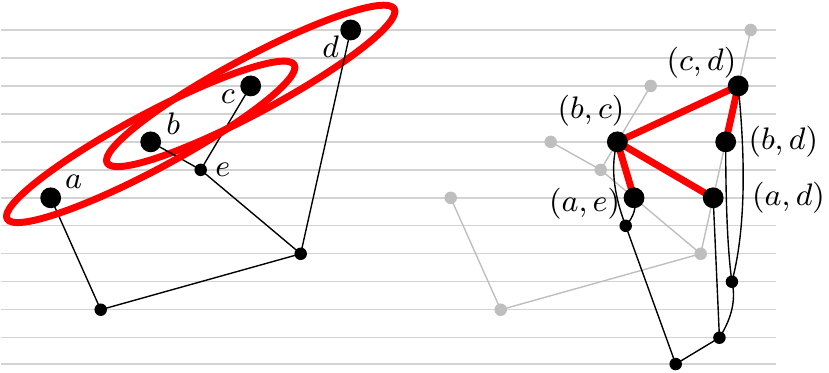}
\caption{Examples of structures determining the big Ramsey degree of a hypergraph on $4$ vertices. The tree $T(V^0,\leq^0,\preceq\restriction_{V^0})$ is on the left and $T(V^1,\leq^1,\preceq\restriction_{V^1})$ with $G^1$ on the right, see Definition~\ref{def51}.}
\label{fig:3}
\end{figure}

Clearly, every $n$-tuple of vertices of $\mathcal R$ has a shape in tree $T^\preceq_\mathcal R$ as well as shapes in trees $T^\preceq_{R_b}$ for every vertex $b$ in the $n$-tuple. All those shapes together determine a very rich (but finite) colouring.
We describe the big Ramsey degree in terms of isomorphism types of structures which combine those trees into a unified object consisting of two trees: first tree is a subtree of $T^\preceq_\mathcal R$ and the second combines subtrees of all relevant trees $T^\preceq_{R_a}$, see also Figure~\ref{fig:3}. Since all the trees are defined by the well-order $\preceq$ we arrive to a perhaps surprising situation where several vertices might be assigned to each vertex $a$ of $\mathcal R$ and thus $\preceq$ turns to a well-preorder.

We are now ready to describe the structures giving the big Ramsey degrees of the 3-uniform hypergraph.
First we describe the tree originating from  $T^\preceq_{\mathcal R}$.

\begin{definition}
Let $(\leq,\preceq)$ be a pair of compatible orders of a set $V'$, let $V$ be the set of leaf vertices of $T(V',\leq,\preceq)$,
and let $\mathcal G=(V,\mathcal E)$ be a $3$-uniform hypergraph.
We say that $\mathcal G$ is \emph{compatible} with $T(V',\leq,\preceq)$ if for every 4-tuple $a,b,c,d$ of distinct vertices of $V$ satisfying $d\preceq c\preceq (a\wedge b)$ we have $\{a,c,d\}\in \mathcal E$ if and only if $\{b,c,d\}\in \mathcal E$.
\end{definition}

Next we introduce a way to synchronize this tree with the union of trees originating from the Rado graphs $R_a$.

\medskip

Given a tree $T(V^0,\leq,\preceq)$ and a compatible 3-uniform hypergraph $\mathcal G=(V,\mathcal E)$, we define
the \emph{neighbourhood graph} of $\mathcal G$ with respect to $T(V^0,\leq,\preceq)$ as the graph $G^1=(V'',E^1)$
constructed as follows:
\begin{enumerate}
 \item $V''$ consists of all pairs $(a,b)$ such that $a\in V$ (by compatibility $V\subseteq V^0$) and $b\in V^0$,
$a\prec b$ and there is no $c\in V^0$, $c\sqsubset b$ such that $a\prec c\prec b$.
 \item $\{(a,b),(c,d)\}\in E^1$ for $a\preceq c$ iff there exists $e\sqsupseteq d$ such that $\{a,c,e\}\in\mathcal E$.
(This is well defined because of the compatibility of $T(V^0,\leq,\preceq)$ and $\mathcal G$.)

\end{enumerate}
For $(a,b)\in V'$, we define its \emph{projection} $\pi:V\times V^0\to V$ by putting $\pi((a,b))=a$. To simplify the notation, we also define $\pi(a)=a$ for every $a\in V$.

\begin{definition}
\label{def51}
The tuple $(V^0,V^1,\preceq,\leq^0,\leq^1)$  is \emph{compatible} with the 3-uniform hypergraph $\mathcal G=(V,\mathcal E)$  iff:
\begin{enumerate}
 \item $V^0\cap V^1=\emptyset$,
 \item $(\leq^0,\preceq\restriction_{V^0})$ is a compatible pair of orders of $V^0$ and
 $T(V^0,\leq^0,\preceq\restriction_{V^0})$ is compatible with $\mathcal G$,
 \item $(\leq^1,\preceq\restriction_{V^1})$ is a compatible pair of orders of $V^1$ and $T(V^1,\leq^1,\preceq\restriction_{V^1})$ is compatible with the neighbourhood graph $G^1=(V^1,E^1)$ of $\mathcal G$ with respect to $T(V^0,\leq^0,\preceq\restriction_{V^0})$,
 \item $\preceq$ is a well pre-order which satisfies $a\neq b$, $a\preceq b$, $b\preceq a\Rightarrow \pi(a)=\pi(b)$, and both projections are defined.
Moreover, whenever $\pi(a)$ and $\pi(b)$ are defined, $\pi(a)\preceq \pi(b)\Rightarrow a\preceq b$. Finally, for $(a,b),(c,d)\in V^1$, we have $((a,b)\wedge (c,d))\prec (b\wedge d)$.
\end{enumerate}
\end{definition}

\begin{theorem}\label{thm:nas}
The big Ramsey degree of a 3-uniform hypergraph $\mathcal G=([n],\allowbreak \mathcal E)$ in $\mathcal R$ is
the number of non-isomorphic structures $([2n-1]\allowbreak \cup V^1,\preceq,\leq^0,\leq^1,\mathcal E,\mathcal P)$ such that
$([2n-1],V^1,\preceq,\leq^0,\leq^1)$ is compatible with $\mathcal E$, $\preceq\restriction_{[2n-1]}$ is a linear order
and $\mathcal P$ consists of all triples $\{a,b,(a,b)\}$ such that $(a,b)$ is a vertex of the neighbourhood graph $G^1$.
\end{theorem}
Proofs of Theorems~\ref{thm:devlin} and~\ref{thm:sauer} use the Milliken tree theorem for
infinite binary branching trees.
The proof of Theorem~\ref{thm:nas} is also based on an application of the
Milliken tree theorem for product of trees. To our best knowledge, this is the first time this product form with
two trees is used. It also seems to be the first combinatorial application of the Milliken
tree theorem for trees with unbounded branching which appears naturally in the tree of types
of 3-uniform hypergraphs (cf~\cite{dodos2016}).








%

\end{document}